\documentclass[11pt,a4paper,notitlepage]{article}
\usepackage[psamsfonts]{amssymb}
\usepackage{amssymb}
\usepackage{geometry} 
\usepackage{float}
\usepackage{booktabs}
\usepackage{algorithm,algorithmic}
\usepackage{amsmath}
\usepackage{pdflscape}
\usepackage{url}
\usepackage{natbib}
\usepackage{graphicx}
\usepackage[utf8]{inputenc}
\allowdisplaybreaks
\begin{document}

\title{A branch-and-Benders-cut algorithm for a bi-objective stochastic facility location problem}
\author{
Sophie N. Parragh$^{1}$ \: Fabien Tricoire$^{1}$ \: Walter J. Gutjahr$^2$  \\[1ex]
 \small $^1$Institute of Production and Logistics Management, Johannes Kepler University Linz\\
 \small Altenberger Straße 69, 4040 Linz, Austria\\
  \small \texttt{\{sophie.parragh,fabien.tricoire\}@jku.at}\\[2ex]
\small $^2$Department of Statistics
and Operations Research, University of Vienna\\
 \small  Oskar-Morgenstern-Platz 1, 1090 Vienna, Austria\\
  \small \texttt{walter.gutjahr@univie.ac.at}
}
\date{}
\maketitle

\begin{abstract}
In many real-world optimization problems, more than one objective plays a role and input parameters are subject to uncertainty. In this paper, motivated by applications in disaster relief and public facility location, we model and solve a bi-objective stochastic facility location problem. The considered objectives are cost and uncovered demand, whereas the demand at the different population centers is uncertain but its probability distribution is known. The latter information is used to produce a set of scenarios. In order to solve the underlying optimization problem, we apply a Benders' type decomposition approach which is known as the L-shaped method for stochastic programming and we embed it into a recently developed branch-and-bound framework for bi-objective integer optimization. We analyze and compare different cut generation schemes and we show how they affect lower bound set computations, so as to identify the best performing approach. Finally, we compare the branch-and-Benders-cut approach to a straight-forward branch-and-bound implementation based on the deterministic equivalent formulation.
\end{abstract}

\section{Introduction}
Facility location problems play an important role in long-term public
infrastructure planning. Prominent examples concern the location of
fire departments, schools, post offices, or hospitals. They are not
only relevant in public (or former public) infrastructure planning
decisions in “regular” planning situations: they are also of concern
in the context of emergency planning, e.g., relief goods distribution
in the aftermath of a disaster or preparation for slow onset disasters
such as droughts. In many of these contexts, accurate demand figures
are not available; assumed demand values rely on estimates, while
their actual realizations depend, e.g., on the severity of the
slow onset disaster, the demographic population development in an
urban district,
etc. Since facility location decisions are usually long-term
investments, the uncertainty involved in the demand figures should
already be taken into account at the planning stage.

Another important issue is that facility location problems often involve several objectives. On the one hand, client-oriented objectives should be optimized. For example, in cases where it is not possible to satisfy the demand to 100 percent, the total covered demand should be as high as possible. On the other hand, cost considerations also play a role. This implies that decision makers face a trade-off between client-oriented and cost-oriented goals. Instead of combining these two usually conflicting measures into one objective function, it is advisable to elucidate their trade-off relationship. Such an approach provides valuable information to the involved stakeholders and allows for better informed decisions. Following this line of thought, in this paper, we model a bi-objective stochastic facility location problem that considers cost and coverage as two competing but concurrently analyzed objectives. Furthermore, we incorporate stochastic information on possible realizations of the considered demand figures in the form of scenarios sampled from probability distributions.

Motivated by recent advances in exact methods for multi-objective integer programming, we solve this problem by combining a recently developed bi-objective branch-and-bound algorithm \citep{parragh2015branch} with the L-shaped method \citep{van1969shaped}, which applies Benders decomposition \citep{benders1962partitioning} to two-stage stochastic programming problems. We integrate several enhancements, such as partial decomposition, and we compare the resulting approach to using a deterministic equivalent formulation within the same branch-and-bound framework. 

This paper is organized as follows.
In Section \ref{sec:relWork}, we give a short overview of related work in the field of bi-objective (stochastic) facility location. In Section \ref{sec:problem}, we define the bi-objective stochastic facility location problem (BOSFLP) that is subject to investigation in this paper and we discuss L-shaped based decomposition approaches of the proposed model. In Section \ref{sec:methods} we explain how we integrate the proposed decomposition schemes into 
the bi-objective branch-and-bound framework. A computational study comparing the different approaches is reported in Section \ref{sec:results}, and Section \ref{sec:conclusion} concludes the paper and provides directions for future research.

\section{Related work}
\label{sec:relWork}

Our problem is a stochastic extension of a bi-objective {\em maximal covering location problem} (MCLP). The MCLP has been introduced in \citet{church1974maximal}.
It consists in finding locations for a set of $p$~facilities in such a way that a maximum population can be served within a pre-defined service distance. Obviously, in this classical formulation, the number of facilities that are to be opened is a cost measure, and the total population that can be served is a measure of demand coverage. Thus, cost occurs in a constraint, whereas the objective represents covered demand. It is natural to extend this problem to a {\em bi-objective} covering location problem (CLP) where both cost (to be minimized) and covered demand (to be maximized) are objectives. Indeed, bi-objective CLPs of this kind have been studied in several papers, see e.g.~\citet{bhaskaran1990multiobjective}, \citet{harewood2002emergency}, \citet{villegas2006solution}, or \citet{Gutjahr20161}; for further articles, we refer to \citet{farahani2010multiple}.

Another strand of literature relevant in the present context addresses bi-objective {\em covering tour problems} (CTPs). One of the oldest CTP models, that of the maximal covering tour problem (MCTP) introduced by \citet{current1994median}, has already been cast in the form of a bi-objective problem. A fixed number of nodes have to be selected out of the
nodes of a given transportation network for being visited by a vehicle, and a tour on these visited
nodes has to be determined. The objectives are minimization of the
total tour length (a cost measure) and maximization of the total demand that is
satisfied within some pre-specified distance from a visited node (a measure of demand coverage).

Other multi-objective CTP formulations can be found in the following papers: \citet{jozefowiez2007bi} deal with a bi-objective CTP where
the second objective function of the MCTP is replaced by the largest
distance between a node of some given set and the nearest
visited node. \citet{doerner2007multicriteria} develop a three-objective CTP model for mobile health
care units in a developing country. \citet{nolz2010bi} study a multi-objective CTP
addressing the problem of delivery of drinking water to
the affected population in a post-disaster situation.

Most importantly for the present work, \citet{tricoire2012bi} generalize the bi-objective CTP to the stochastic case by assuming uncertainty on demand. The aim is to support the choice of distribution centers (DCs) for relief commodities and of delivery tours supplying the DCs from a central depot. Demand in population nodes is assumed as uncertain and modeled stochastically. DCs have fixed given capacities, as well as vehicles. The model considers two objective functions: The first objective is cost (more precisely, the sum of opening costs for DCs and of transportation costs), and the second is expected uncovered demand. Contrary to basic covering tour models supposing a fixed distance threshold, uncovered demand is defined by the more general model assumption that the percentage of individuals who are able and willing to go to the nearest DC can be represented by a nonincreasing function of the distance to this DC. Both the demand of those individuals who stay at home and the demand of those individuals who are not supplied in a DC because of DC capacity and/or vehicle capacity limits contribute to the total uncovered demand. Because of the uncertainty on the actual demand, total uncovered demand is a random variable, the expected value of which defines the second objective function to be minimized.

The problem investigated in the present paper generalizes the bi-objective CLP
to a stochastic bi-objective problem by considering demand as uncertain and
modelling it by a probability distribution, in an analogous way as in \citet{tricoire2012bi}. Alternatively, the investigated problem can also be derived from a CTP by omitting the routing decisions and generalizing the resulting bi-objective location problem again to the stochastic case. From the viewpoint of the latter consideration, the current model can be seen as related to the special case of the problem of \citet{tricoire2012bi} obtained by neglecting routing costs. However, the current model builds on refined assumptions concerning the decision structure of the two-stage stochastic program which makes the second-stage optimization problem nontrivial, contrary to \citet{tricoire2012bi} where the second-stage optimization problem can be solved by elementary calculations.

Multi-objective stochastic optimization (MOSO) problems, though of eminent importance for diverse practical applications, are investigated in a more limited number of publications, compared to the vast amount of literature both on multi-objective optimization and on stochastic optimization; for surveys on MOSO, we refer the reader to \citet{CCMR04}, \citet{Abdelaziz20121}, and \citet{gutjahr2013stochastic}. Of special relevance for our present work are multi-objective {\em two-stage stochastic programming} models where the ``multicriteria'' solution concept is that of the determination of {\em Pareto-efficient solutions}, and where the first-stage decision contains {\em integer} decision variables. Most papers in this area assume that one of the two objectives only depends on the first-stage decision, whereas the other objective depends on the decisions in both stages. This holds also for \citet{tricoire2012bi}. Let us give two other examples: \citet{FGOS10} present a two-stage stochastic bi-objective mixed integer program for reverse logistics, with strategic decisions on location and function of diverse collection and recovery centers in the first stage, and tactical decisions on the flow of disposal from clients to centers or between centers in the second stage. The first objective is expected total cost, which depends both on the first-stage and second-stage decision, whereas the second objective, the obnoxious effect on the environment, only depends on the first-stage decision. Stochasticity is associated with waste generation and with transportation costs. \citet{CVAO11} deal with decisions on the location of DCs and on transportation flows in a two-echelon production distribution network. Uncertainty holds with respect to the demand. The first objective represents expected costs, whereas the second objective expresses the sum of the maximum lead times from plants to customers. The authors model the random distribution by scenarios and solve the two-stage programming model by the L-shaped method, a technique that we will also use in our present work.

Our problem has a structure similar to the models cited above: while the covered demand depends on the decisions in both stages, the cost objective is already determined by the first-stage decision. This allows the development of an efficient solution algorithm. Contrary to \citet{tricoire2012bi}, we will not apply an epsilon-constraint method for the determination of Pareto-efficient solutions, but use instead a more recent method developed in \citet{parragh2015branch}.

Finally, let us mention that the humanitarian logistics literature,
which tackles facility location problems under high
uncertainty and multiple objectives, has been relatively prolific
with regards to multi-objective stochastic optimization models and
corresponding solution techniques (cf.~\cite{GutjahrNolz2015}). Let us
give two examples of papers using both Pareto optimization and
two-stage stochastic programs. \citet{khorsi2013nonlinear} propose a
bi-objective model with objectives ``weighted unsatisfied demand'' and
``expected cost''. Discrete scenarios from a set of possible disaster
situations are applied to represent uncertainty. The
epsilon-constraint method is used to solve the
model. \citet{rath2015bi} deal with the uncertain accessibility of
transportation links and develop a two-stage stochastic programming
model where in the first stage decisions on the locations of
distribution centers have to be made, and in the second stage (based
on current information on road availability) the transportation flows
have to be organized. Objective functions are expected total cost and
expected covered demand. The structure of the two-stage stochastic
program is different from that in the current paper insofar as both
objectives depend on both first-stage and second-stage decision
variables, which requires specific (and computationally less
efficient) solution techniques.
For general information on humanitarian logistics, the reader is
referred to the standard textbook by
\citet{tomasini2009preparedness}. Two-stage stochastic programing
approaches to this field (in a single-objective context) are reviewed
in \citet{grass2016two}. A good recent example for the application
of Benders decomposition to a stochastic model for humanitarian relief
network design is \citet{elcci2018chance}.

For an overview on facility location in general, we refer the reader to \citet{hamacher2002facility}. Standard textbooks on multi-objective optimization and on stochastic programming are \citet{Ehrg05} and \citet{birge2011introduction}, respectively.

\section{Problem definition and decomposition}
\label{sec:problem}
In the bi-objective stochastic facility location problem (BOSFLP)
considered in this article, the demand at each node $i \in V$ is uncertain. We denote by $W_i$ the demand at node $i$, by ${\mathbf E}(W_i)=w_i$ the expected demand at node $i$ and  by $\xi_i$ a  random variable such that $W_i = \xi_i w_i$. At each node $j$, a facility may be built. A facility at node $j$ has a capacity $\gamma_j$ and operating costs $c_j$. Furthermore, facilities that are farther than a certain maximum distance $d_{max}$ from a demand point may not be used to cover it. In order to take this aspect into account, we consider the set $A = \{(i,j)| i,j \in V, d_{ij} \leq d_{max} \}$ of possible assignments $(i,j)$, where $d_{ij}$ denotes the distance of demand node $i$ from a potential facility at node $j$. The two considered goals are to minimize the total costs for operating facilities and to maximize the expected covered demand. 
Using the following decision variables,
%
\begin{align*}
	z_j  = &\begin{cases}
			1$ if a facility is built at node $j$ and $
			\cr
			0$, otherwise.$
	\end{cases}\\
	y_{ij} &\text{ demand of population node $i$ covered by facility $j$},\\
	u_j &\text{ total demand covered by facility $j$.}
\end{align*}
we formulate the BOSFLP as a two-stage stochastic program:\\
\begin{align}
	\min f_1 & = \sum_{j \in V} c_j z_{j} 
		& \label{mod1:of1}\\
	\min f_2 & = {\mathbf E}( Q(z,\xi) ) 
		&\label{mod1:of2}
\end{align}
\begin{align}
	z_{j} & \in \{0,1\} && \forall j \in V.
\end{align}
Second stage:
\begin{align}
	Q(z,\xi) & = \min_{u} ( - \sum_{j \in V} u_j)
	\label{mod2:of}
\end{align}
subject to:
\begin{align}
	u_j & \leq \sum_{i | (i,j) \in A} y_{ij}
		&& \forall j \in V,
		\label{mod2:cov}\\
	u_j & \leq \gamma_j z_j
		&& \forall j \in V,
		\label{mod2:cap}\\
	y_{ij} & \leq W_i z_{j} 
		&& \forall (i,j) \in A,
		\label{mod2:open}\\
	\sum_{j | (i,j) \in A} y_{ij} & \leq W_i  
		&& \forall i \in V, 
		\label{mod2:ySumUpToOne}\\
	0 & \leq y_{ij} 
		&& \forall (i,j) \in A,\\
	0 & \leq u_j 
		&& \forall j \in V.
\end{align}

Objective function \eqref{mod1:of1} minimizes the total facility
opening costs. Objective function \eqref{mod1:of2} maximizes the
expected covered demand. In order to obtain two minimization
objectives, objective \eqref{mod1:of2} has been multiplied by $(-1)$.
Note that maximization of the expected covered demand is equivalent to
minimization of the expected uncovered demand, since the expected
demand is a constant.
The first stage model only comprises one set of constraints which
require that all $z_j$ variables may only take values $0$ or $1$. The
second stage model consists of the objective function given in
\eqref{mod2:of}, representing the negative value of the total covered
demand, and a number of constraints which determine the maximum
possible coverage given a first stage solution. Constraints
\eqref{mod2:cov} link the coverage variables with the assignment
variables: the covered demand at node $j$ cannot be larger than the
actual demand assigned to this node. Constraints \eqref{mod2:cap} make
sure that the capacity of facility $j$ is not exceeded. Constraints
\eqref{mod2:open} guarantee that a demand node can only be assigned to
a facility if the respective facility is open. Finally, constraints
\eqref{mod2:ySumUpToOne} ensure that any part of the demand at $i$ is
only covered at most once. The variables $z_j$ are first-stage
decision variables whereas the variables $y_{ij}$ are second-stage
decision variables, i.e. the latter variables can depend on the
realizations of the demand values; that is, $y_i=y_i(\xi)$. Similarly,
$u_j=u_j(\xi)$ and $W_i=W_i(\xi)$.

In this paper, we use a discrete set of scenarios (with equal probabilities) in order to approximate the (joint) probability distribution of the demand as estimated by the decision maker. If a Monte-Carlo simulation procedure is available for simulating demand, each realization of this procedure can be taken as a scenario. Let $N$ denote this set of scenarios. Then, using an additional index $\nu$ to denote a given scenario for the variables of the second stage problem, we obtain the following expanded or deterministic equivalent model:

\begin{align}
	\min f_1 & = \sum_{j \in V} c_j z_{j} 
		& \label{mod2:of1-2}\\
	\min f_2 & = - \frac{1}{|N|} \sum_{\nu \in N} \sum_{j \in V} u_j^{\nu}
		&\label{mod2:of2-2}
\end{align}
subject to:
\begin{align}
	u_j^{\nu} & \leq \sum_{i | (i,j) \in A} y_{ij}^{\nu}
		&& \forall j \in V, \nu \in N,
		\label{mod2:cov-2}\\
	u_j^{\nu} & \leq \gamma_j z_j
		&& \forall j \in V, \nu \in N,
		\label{mod2:cap-2}\\
	y_{ij}^{\nu} & \leq W_i ^{\nu}  z_{j}
		&& \forall (i,j) \in A, \nu \in N,
		\label{mod2:open-2}\\
	\sum_{j | (i,j)\in A} y_{ij}^{\nu} & \leq W_i ^{\nu} 
		&& \forall i  \in V, {\nu} \in N, 
		\label{mod2:ySumUpToOne-2}\\
	0 & \leq y_{ij}^{\nu} 
		&& \forall (i,j) \in A, \nu \in N, \\
	0 & \leq u_j^{\nu} 
		&& \forall j \in V, \nu \in N ,\\
	z_{j} & \in \{0,1\} && \forall j \in V.
\end{align}

This problem formulation decomposes by scenario and the L-shaped method, as introduced by \citet{van1969shaped}, can be used to solve its linear relaxation. The L-shaped method relies on a master problem and one subproblem per scenario whereas the information from the subproblem is incorporated into the master problem by means of cutting planes. We distinguish feasibility and optimality cuts. In the case of complete recourse, as is the case for our problem, only optimality cuts have to be added: the solution to the first stage problem will always allow a feasible solution to the second stage problem. (This holds for our problem, since setting $y_{ij}$ and $u_j$ to zero produces a feasible solution.)

More precisely, using variable $\theta$ to represent the second stage objective, 
we obtain the following master linear program (LP):

\begin{align}
	\min f_1 = &\sum_{j \in V} c_j z_j\\
	\min f_2 = &\: \theta
		\label{LSM:of}
\end{align}
subject to:
\begin{align}
	- \sum_{i \in V} \max_{\nu \in N} \{ W^{\nu}_i \}& \leq \theta \\ 
  0  \leq z_j  & \leq 1  && \forall j \in V 
  	\label{LSM:nonnegTheta}
\end{align}
where $- \sum\limits_{i \in V} \max\limits_{\nu \in N} \{ W^{\nu}_i \}$ provides a valid bound on $\theta$, since $Q(z,\xi)  = \min\limits_{u} ( - \sum\limits_{j \in V} u_j) \geq - \sum\limits_{i \in V} \max\limits_{\nu \in N} \{ W^{\nu}_i \}$.

As explained in Section~\ref{sec:method}, the two objectives are
combined into a weighted sum and the resulting single-objective LP is
iteratively solved with different weights in order to enumerate the
set of efficient solutions. In that context and for a given set of
weights, to determine if the obtained solution to the master
weighted-sum LP is optimal, we check if optimality cuts have to be
added. We denote by $z_j^l$ and $\theta^l$ the
variable values obtained from solving the master LP and we solve for each
$\nu \in N$ the following model:
\begin{align}
	\min (- \sum_{j \in V} u_j)
		\label{LSstep3:of}
\end{align}
subject to:
\begin{align}
	u_j & \leq \sum_{i | (i,j) \in A} y_{ij}
		&& \forall j \in V,
		\label{LSstep3:distanceConstr}\\
	u_j & \leq \gamma_j z^l_j
		&& \forall j \in V,
		\label{LSstep3:capacityConstr}\\
	y_{ij} & \leq W_{i}^{\nu} z^l_{j} 
		&& \forall (i,j) \in A, 
		\label{LSstep3:openingConstr}\\
	\sum_{j | (i,j) \in A} y_{ij} & \leq	W_{i}^{\nu}
		&& \forall i \in V, 
		\label{LSstep3:populationConstr}\\
	0 & \leq y_{ij} 
		&& \forall (i,j) \in A,
		\label{LSstep3:yBounds}\\
	0 &\leq u_j
		&& \forall j \in V.
		\label{LSstep3:uBounds}
\end{align}
Let $\mathcal{Q}(z) = {\mathbf E}_{\xi}Q(z,\xi)$. If $\mathcal{Q}(z^l) \leq
\theta^l$, we terminate: optimality has been reached. Otherwise, we generate an optimality cut. Optimality cuts rely on dual information. To write the dual of the above model \eqref{LSstep3:of} -- \eqref{LSstep3:uBounds}, we denote by $\lambda_j$ the dual variables of constraints \eqref{LSstep3:distanceConstr}, by $\pi^l_j$ the dual variables
of constraints \eqref{LSstep3:capacityConstr}, by $\sigma_{ij}$ the dual variables of constraints
\eqref{LSstep3:openingConstr}, and by $\delta_i$ the dual variables of constraints \eqref{LSstep3:populationConstr}:
\begin{align}
	\max (- \sum_{j \in V} \pi_j \gamma_j z_j^l -  \sum_{(i,j) \in A} \sigma_{ij} W_{i}^{\nu} z_j^l - \sum_{i \in V} W_{i}^{\nu} \delta_i)
		\label{dual:of}
\end{align}
subject to:
\begin{align}
	\lambda_j + \pi_j  & \geq 1
		&& \forall j \in V,
		\label{dual:Constr1}\\
	\lambda_j - \sigma_{ij} - \delta_i    & \leq 0
		&& \forall (i,j) \in A,
		\label{dual:Constr2}\\
	\lambda_j &\geq 0
		&& \forall j \in V,
		\label{dual:lBoundsL}\\
	\pi_j &\geq 0
		&& \forall j \in V,
		\label{dual:lBoundsP}\\
	\sigma_{ij} &\geq 0
		&& \forall (i,j) \in A,
		\label{dual:lBoundsS}\\ 
	\delta_i & \geq 0
		&& \forall i \in V,
		\label{dual:lBoundsD}
\end{align}
With a bit of abuse of notation, we denote by $\pi_j(\xi^{\nu})$,  $\sigma_{ij}(\xi^{\nu})$, and $\delta_i(\xi^{\nu})$  the dual variable values for a given scenario $\nu$ and the
objective function value by 
$Q(z^l, \xi^{\nu}) = - \sum_{j \in V} \pi_j(\xi^{\nu}) \gamma_j z_j^l -  \sum_{(i,j) \in A} \sigma_{ij}(\xi^{\nu}) W_{i}^{\nu} z_j^l - \sum_{i \in V} W_{i}^{\nu} \delta_i(\xi^{\nu})$. 
Then, the optimality cut is of the following form:
\begin{align}
	Q(z,\xi^{\nu})  \geq  Q(z^l,\xi^{\nu})
	&+\big( \big[ - \sum_{j \in V} \pi_j(\xi^{\nu}) \gamma_j z_j - 
	 \sum_{(i,j) \in A} \sigma_{ij}(\xi^{\nu}) W_{i}^{\nu} z_j  - 
	\sum_{i \in V}  \delta_i(\xi^{\nu}) W_{i}^{\nu} \big] \\ \nonumber
	& - \big[ - \sum_{j \in V} \pi_j(\xi^{\nu}) \gamma_j z_j^l -  
	\sum_{(i,j) \in A} \sigma_{ij}(\xi^{\nu}) W_{i}^{\nu} z_j^l 
	- \sum_{i \in V}  \delta_i(\xi^{\nu}) W_{i}^{\nu} \big]  \big)
\end{align}
Rearranging the terms, we obtain
\begin{align}
	Q(z,\xi^{\nu})  \geq  Q(z^l,\xi^{\nu})
	 + \big( \sum_{j \in V}
  \pi_j(\xi^{\nu}) \gamma_j (z_j^l - z_j)  + \sum_{(i,j) \in A}
  \sigma_{ij}(\xi^{\nu}) W_{i}^{\nu} (z_j^l - z_j) \big)
  \label{LScut:indiv}
\end{align}
%
Finally, combining over all scenarios, we obtain the optimality cut for the expected
second stage objective function (with $\bar{\pi}_j = 1/|N| \sum_{\nu}
\pi_j(\xi^{\nu})$ and 
$\bar{\sigma}_{ij} =1/|N| \sum_{\nu} \sigma_{ij}
(\xi^{\nu}) W_{i}^{\nu}$):
\begin{align}
  \theta \geq   \mathcal{Q}(z^l) + \big( \sum_{j \in V} \bar{\pi}_j \gamma_j
  (z_j^l - z_j)  + \sum_{(i,j) \in A} \bar{\sigma}_{ij} (z_j^l - z_j) \big)
  \label{LScut:avg}
\end{align}

Alternatively, instead of adding one optimality cut per iteration, we can also add one cut per scenario. In order to do so, a separate variable $\theta^{\nu}$ for each scenario $\nu$ has to be used, resulting in the following master LP:
\begin{align}
	\min f_1 = & \sum_{j \in V} c_j z_j \\
	\min f_2 = & \frac{1}{N} \sum_{\nu \in N}  \theta^{\nu}
		\label{LSMindiv:of}
\end{align}
subject to:
\begin{align}
  - \sum_{i \in V} W_i^{\nu} & \leq \theta^{\nu} \leq 0
  		&& \forall \nu \in N
  	\label{LSMindiv:nonnegTheta}\\
  	 0 & \leq z_j \leq 1  
  		&& \forall j \in V
\end{align}
Then, we check each subproblem and add a cut of the form \eqref{LScut:indiv} in the case where $\theta^{\nu,l} < Q(z^l,\xi^{\nu})$. In the case where all scenarios are checked and no additional cut has to be added, optimality has been reached.
%

\section{Solution methods}
\label{sec:methods}
\label{sec:method}

In order to solve the BOSFLP, we integrate L-shaped based cut generation into the recently introduced bi-objective branch-and-bound framework of \citet{parragh2015branch}. In what follows, we first describe the key ingredients of the branch-and-bound framework and thereafter how we combine it with the L-shaped method.

\subsection{Bi-objective branch-and-bound}
Without loss of generality, we consider a bi-objective
\emph{minimization} problem.
The bi-objective branch-and-bound (BIOBAB)
algorithm of \citet{parragh2015branch} generalizes the
single-objective concept of branch-and-bound to two objectives. We
first introduce the notion of lower and upper bound set. Thereafter,
we explain the main loop of the algorithm and we describe several
enhancements.

\subsubsection{Lower and upper bound sets}
\label{sec:boundsets}
During the execution of the BIOBAB algorithm, instead of single numerical values, upper and lower bound sets are
computed. These rely on the notion of \emph{bound set} introduced by
\citet{Ehrgott2006}. A subset $L$ of the objective space is called a
lower bound set of the feasible set $Y$ in objective space if $\forall
y \in Y \exists x \in L: y \geq x$, where $y \geq x$ iff $y_i \geq x_i
(i=1, 2)$.
Starting from the root node, at each node of the branch-and-bound tree, a lower bound (LB) set is calculated. The 
special LB set we use corresponds to the lower left boundary of the
convex hull of the feasible set of the current node LP in objective
space. This boundary can be described by its corner points which can
be efficiently computed by means of an algorithm that is similar to
that of~\cite{Aneja1979}.
This algorithm consists in solving a series of single-objective
weighted-sum problems by systematically enumerating a finite set of
weight combinations. In our case, this means that the two objectives
are combined into a weighted sum and we solve a (linear) relaxation of
the weighted-sum problem with the appropriate weight combination at
every step of the LB generation algorithm. The image of the solution
to each relaxed weighted-sum problem gives a corner point of the
boundary in objective space. In a first step, the algorithm computes
the two extreme solutions, i.e., the best solution optimizing $f_1$
and the best solution optimizing $f_2$. Let $a$ denote the point in
objective space which is the image of the optimal solution for $f_1$
and $b$ the point in objective space which is the image of the optimal
solution for $f_2$. In order to obtain the best possible value for the
respective other objective function, lexicographic minimization is
used. The next step consists in identifying the weights for finding
the next solution of the LB set. These weights are derived from $a$
and $b$. Let $a_1$ and $a_2$ denote the coordinates of $a$ and $b_1$
and $b_2$ the coordinates of $b$. Then the weights to obtain the next
solution of the LB set are $w_2 = b_1 - a_1$ and $w_1 = a_2 -
b_2$. Let the image of this new solution in objective space be denoted
by $c$. Then we look for additional solutions between $a$ and $c$ and
between $c$ and $b$, in the same way as before. For further details we
refer to \citet{parragh2015branch}.

\begin{figure}
\begin{center}
\includegraphics[width=0.5\textwidth]{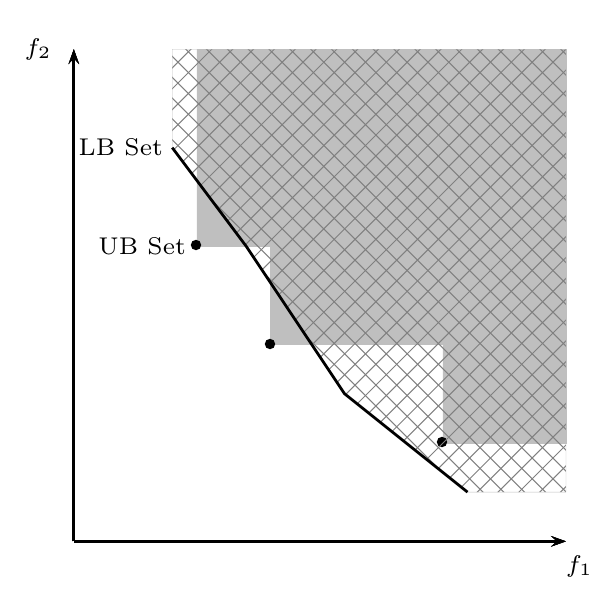}
\end{center}
\caption{LB set (continuous) of a given branch-and-bound node, current
  UB set (dots), region dominated by LB set (grid pattern), region
  dominated by UB set (gray) in objective space. Non-dominated
  feasible solutions have to lie in the white grid-pattern area.}
\label{fig:OSB}
\end{figure}

The thus obtained LB set is then filtered using a set of known
solutions, called the upper bound (UB) set. The UB
set corresponds to all integer feasible solutions obtained during the
search that have not been found to be dominated so far. In order to
fathom a node, the whole LB set of this node must be dominated by the
UB set. Note that the LB set and the UB set are conceptually
different; while the former is a continuous set, the latter is a
discrete set. This is illustrated in Figure \ref{fig:OSB}. In the
example of Figure \ref{fig:OSB}, the current node cannot be fathomed
since the UB set does not dominate the LB set.

\subsubsection{Tree generation and branching rules}
Algorithm~\ref{alg:biobab} shows the main loop of the BIOBAB algorithm. Function $push(C, x)$ adds $x$ to an existing collection of nodes
$C$ and function $pop(C)$ retrieves a node from $C$. They are used to add and retrieve the nodes of the
branch-and-bound tree. 
A node in the branch-and-bound tree represents a set of branching decisions. Depending on the data structure employed for $C$, different tree exploration strategies can be obtained. We use depth-first search. 
The algorithm can take a starting UB set as input in order to speed up the search. In our case, the UB set passed to the algorithm is empty and it is updated every time a new integer solution is found.

\begin{algorithm}
  \caption{$BIOBAB(UB)$}
  \begin{algorithmic}[1]
    \STATE $rootNode \gets \emptyset$
    \STATE $C \gets \emptyset$
    \STATE $push(C, rootNode)$
    \WHILE{$C \neq \emptyset$}
    \STATE $node \gets pop(C)$
    \STATE $LB \gets bound(node, UB)$
    \STATE $LB \gets filter(LB,UB)$ \label{line:filtering}
    \IF {$LB \neq \emptyset$}
      \STATE $newBranches \gets branch(LB)$
      \FORALL{$decision \in newBranches$}
      \STATE $push(C, node \cup \{decision\})$
      \ENDFOR
    \ENDIF
    \ENDWHILE
    \RETURN $UB$
  \end{algorithmic}
  \label{alg:biobab}
\end{algorithm}

A key component of the BIOBAB algorithm of \citet{parragh2015branch} is a branching rule that works on the objective space, which is referred to as \emph{objective space
branching}. It allows to discard dominated regions of the search space even if a given node cannot be fathomed. The information whether or not objective space branching can be performed is obtained in the filtering step: whenever the current UB set allows to discard regions from the lower bound set and the resulting LB set is discontinuous, objective space branching in combination with variable branching is performed and each new branch (or node) corresponds to a different continuous subset of the discontinuous LB set. In the example depicted in Figure \ref{fig:OSB}, the filtering operation results in a discontinuous LB set consisting of three continuous subsets.

Whenever the filtered LB set is not discontinuous, standard variable
branching rules are applied. In this case, the binary variable to
branch on is selected based on information from all corner point
solutions of the LB set: the variable that is fractional in the
highest amount of corner points is selected for branching. Ties are
broken by selecting the variable with lowest average distance to
0.5. If this is not enough, ties are broken by selecting the variable
whose average value has the lowest distance to 0.5.

\subsubsection{Enhanced objective space filtering}
Another key component of the algorithm of \citet{parragh2015branch} are enhanced objective space filtering rules that rely on the observation that the objective values of integer solutions may only take certain values. In the simplest case, they are restricted to integer values. In \citet{parragh2015branch}, only integer problems are addressed, where all coefficients in the objective function may only assume integer values. In this case, it is easy to observe that integer solutions may only assume integer objective values. In this paper, we solve a mixed integer program (the $y_{ij}$ and $u_j$ variables may assume fractional values). However, the continuous variables only appear in the second objective function. Thus, for the first objective function the same reasoning as in  \citet{parragh2015branch} can be used. Our second objective function depends on continuous variables and, in addition, we divide by the number of scenarios to obtain the expected value.  However, we can still exploit the ideas of  \citet{parragh2015branch}. The reasoning is as follows. Let us assume that all coefficients are integer valued (both in the constraints and in the objective functions). This implies that the capacities $\gamma_{j}$ of the distribution centers are integer valued as well as the demands at the demand nodes $W_i$. Then, it is easy to see that, in any optimal solution for a given scenario, at each distribution center, either the capacities are fully used (we maximize covered demand) or, in the case of excessive capacities, the entire demand of the reachable demand nodes is covered, resulting in an integer valued objective function. Now, fractional values can only be due to the term $1/|N|$. Since this term is constant, we can simply multiply the second objective function by $|N|$ to obtain integer valued results. If we do not want to do that, the constant term $1/|N|$ still allows us to know the granularity of the admissible values of the objective function; any region in the objective space which does not contain any admissible values can be removed from further consideration. This observation can be used to prune LB segments and to speed up the LB set computation procedure. For further details we refer to \citet{parragh2015branch}.



\subsection{Lower bound set generation and integration with L-shaped method}


Integrating the L-shaped method into BIOBAB mainly affects the lower bound set generation scheme. In what follows we first present the employed master program and then the employed cut generation strategies.

\subsubsection{Master program}
\label{sec:master}

In a first step, we set up the master LP. 
Let $w_1$ and $w_2$ denote the weights as described
in~\ref{sec:boundsets} and $\bar{f}_1$ and $\bar{f}_2$ the upper
bounds on $f_1$ and $f_2$ respectively, we obtain the following
generic master LP for the multi-cut version:
%
\begin{align}
	\min w_1 \sum_{j \in V} c_j z_j + w_2 \frac{1}{|N|} \sum_{\nu \in N} \theta^{\nu}
\end{align}
subject to:
\begin{align}
  - \sum_{i \in V} W_i^{\nu}  \leq \theta^{\nu} & \leq 0
  		&& \forall \nu \in N
  	\label{LSMindiv:nonnegTheta-2}\\
  	\sum_{j \in V} c_j z_j & \leq \bar{f}_1 \\
  	\frac{1}{|N|} \sum_{\nu \in N}\theta^{\nu} & \leq \bar{f}_2 \\
  	 0 \leq z_j & \leq 1  
  		&& \forall j \in V
\end{align}
To obtain the single cut version of the master LP, $1/|N|
\sum\limits_{\nu \in N} \theta^{\nu}$ has to be replaced by
$\theta$. However, preliminary experiments indicated that, as
expected, the multi-cut version performs better than the single cut
version. For that reason, we focus on the multi-cut version. The model features bounds on both objectives to allow for easy updates in the case of objective space branching, which is realized by updating these bounds to discard dominated regions of the objective space.
The weights in the objective function are determined by the algorithm of \citet{Aneja1979}.  
For each weight combination the L-shaped method is applied,
i.e. optimality cuts (see Section~\ref{sec:problem}) are generated as explained in the subsequent Section. 

In order to strenghten the above master LP, we can use the following valid inequalities:
\begin{align}
 -  \sum_{j \in V} z_j \gamma_j & \leq \theta^{\nu} && \forall \nu \in N
\end{align}
They rely on the fact that the maximum coverage level is bounded by the total capacity of the number of opened facilities. By doing so, in the case where capacities are tight, we anticipate that fewer optimality cuts have to be added.

\subsubsection{Optimality cut generation strategies}

In the general case, the master program is solved, each scenario subproblem is solved and in the case where $\theta^{\nu}$ is currently under-estimated an optimality cut is added and the master program is solved again. Optimality is attained if no additional cut has to be added. However, it is clearly not necessary to check all scenarios for valid cuts at each iteration: we can stop cut generation as soon as at least one cut has been found. In order to do so, several strategies can be envisaged and our preliminary experiments showed that the following strategy works reasonably well:
at each call to the optimality cut generation routine, we do not start
to check for cuts with the first scenario but we start with the
scenario following the last scenario for which a cut was generated,
i.e. if scenario 2 generated the last cut, in the next iteration we
check scenario 3 first and iterate over the scenarios such that
scenario 2 is the one checked last. This way, we always check first
one of the scenarios that have not been checked for the longest time. In terms of cut management, we maintain a global cut pool and we keep all generated cuts in this pool.

\label{sec:incumbent}
In the single objective case it has been observed, e.g. by
\citet{adulyasak2012benders}, that considerable performance gains are
achieved if optimality cuts are only generated at incumbent
solutions. Motivated by the success in the single objective domain, we
transfer this idea to multi-objective optimization. We recall that we
generate bound sets which are obtained by systematically solving a
series of weighted-sum problems. In the current context, each
weighted-sum problem corresponds to solving the above master program
with the L-shaped method. In the case where fewer optimality cuts than
necessary or even no optimality cuts are added, objective two is
under-estimated and therefore we have a valid lower bound on the true
value of objective two. This means that we do not really need to
generate optimality cuts at each weighted-sum solution but we can
restrict cut generation to those weighted-sum solutions which are
integer feasible (mimicking the idea of adding cuts only at incumbent
solutions). In the following we denote such a solution as an incumbent
solution.

\begin{figure}
\begin{center}
\includegraphics[width=0.5\textwidth]{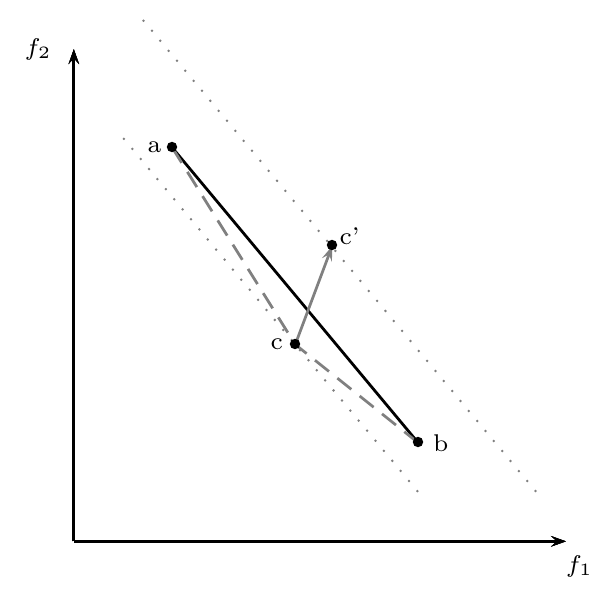}
\end{center}
\caption{LB points $a$ and $b$, new ``incumbent'' point $c$, through cut generation lifted to  $c'$ } 
\label{fig:LBset}
\end{figure}

In the case where we find an incumbent solution, we do
generate cuts then re-solve the modified LP with the same set of
weights, in a loop, until one of two things happens:
\begin{enumerate}
\item the solution is not integer any more
\item the solution is integer and no more cuts can be generated
\end{enumerate}
The point thus obtained is then used as usual for LB set calculation
purposes. 
Figure~\ref{fig:LBset} depicts the situation where points $a$ and $b$ 
have been generated during LB set generation and the next step consists in
investigating the segment between $a$ and $b$. For this purpose the objective 
weights are set to $w_2 = b_1 - a_1$ and $w_1 = a_2 - b_2$ as described above
and we obtain point $c$. Without cutting plane generation, the segments $(a,c)$ and
$(c,b)$ would be investigated (dashed lines in the figure) by the LB set generation scheme. 
Now let us assume
that $c$ is an incumbent. This means that optimality cuts are generated and the cut 
generation loop results in a solution whose image in objective space is the point $c'$.

This point is \emph{above} line $(a,
b)$. 
This is an
issue, as the LB set algorithm only expects 
points below or on that line; this can lead to a non-convex LB
set, because not every point in the convex hull boundary used the same
cuts, i.e. the LP changed during the process. However the branch-and-bound
algorithm relies on a convex LB set. Therefore, in such an
eventuality, the new point is discarded for LB set calculation
purpose, and the segment $(a, b)$ is kept as valid (albeit not tight)
LB segment.
The LB segment is valid since the objective function value level curve of $c'$, depicted by a dotted line in Figure~\ref{fig:LBset}, is a valid LB (set) \citep{stidsen2014branch}.


\subsubsection{Partial decomposition}
Following \citet{crainic2014partial}, partial decomposition appears
to be a viable option to obtain further speedups in the context of a Benders type algorithm. It refers to incorporating some of the scenarios into
the master problem. Let $N^*$ denote the set of scenarios that are incorporated
into the master LP, different strategies regarding which scenario should be
part of $N^*$ can be envisaged. In the simplest case, the first scenario is put
into $N^*$. 
After preliminary testing we decided to keep the scenario with lowest
deviation from the average scenario, plus the $k$ scenarios with
highest deviation.

\section{Computational experiments}
\label{sec:results}
The previously described algorithms have been implemented using Python and
Gurobi 8. The algorithms are run on a cluster with Xeon E5-2650v2 CPUs
at 2.6 GHz. Each job is allocated 8 GB of memory and two hours of CPU
effort. Multi-threading is disabled in Gurobi.
In what follows, we first give an overview of the considered benchmark
instances. Thereafter, we compare the different methods and we discuss
the obtained results.

\subsection{Benchmark instances}

We use a set of 26 instances which are derived from real world data
from the region of Thi\`{e}s in western Senegal (for further details on
these data we refer to \citet{tricoire2012bi}). These instances
feature between 9 and 29 vertices. Only 10 scenarios were used
in~\citet{tricoire2012bi}; we use 10, 50, 100, 200, 300, 400, 500, 600,
700, 800, 900 and 1000 scenarios for each instance. Scenarios are
generated using a procedure similar to the one
in~\citet{tricoire2012bi}. Using more scenarios improves the quality
of the approximation of the real situation but typically requires
additional CPU effort.

\subsection{Cut separation settings}
We compare several settings:
\begin{itemize}
\item \emph{no decomposition}: all constraints are considered
  explicitly in the master problem, no cuts are generated during the
  branch-and-bound algorithm,
\item \emph{base}: base L-shaped method. No partial decomposition, no
  valid inequalities, cuts are systematically generated when they are
  violated.
\item \emph{partial decomposition}: the scenario with lowest deviation
  from average is built in the master problem, as well as the 4
  samples with highest deviation.
\item \emph{valid inequalities}: the valid inequalities described at
  the end of  Section~\ref{sec:master} are added to the master
  problem.
\item \emph{incumbent cuts}: cuts are generated systematically at the
  root node of the tree search, then only on incumbent
  solutions, as described in Section~\ref{sec:incumbent}.
\item \emph{incumbent cuts + valid inequalities}: Both strategies are used.
\end{itemize}

We first compare all settings in terms of CPU effort. Since there are
6 settings, 26 instances and 12 sample sizes, there are 1872 runs to
compare. For that reason we present \emph{performance profiles}. A performance
profile is a chart that compares the performance of various
algorithms~\citep{Dolan2002}. The performance of a setting for a given
instance is the ratio of the CPU effort required with this setting for
this instance over the best known CPU effort for the same
instance. The best performance achievable is always 1. On a
performance profile, performance is indicated on the $x$-axis while
the $y$-axis indicates the ratio of instances solved with at least
that level of performance by a certain setting. If a certain setting
does not converge in solving a given instance within the
allotted CPU budget, then this setting does not provide a performance
for that instance.

\begin{figure}
  \centering
  \includegraphics[scale=.7]{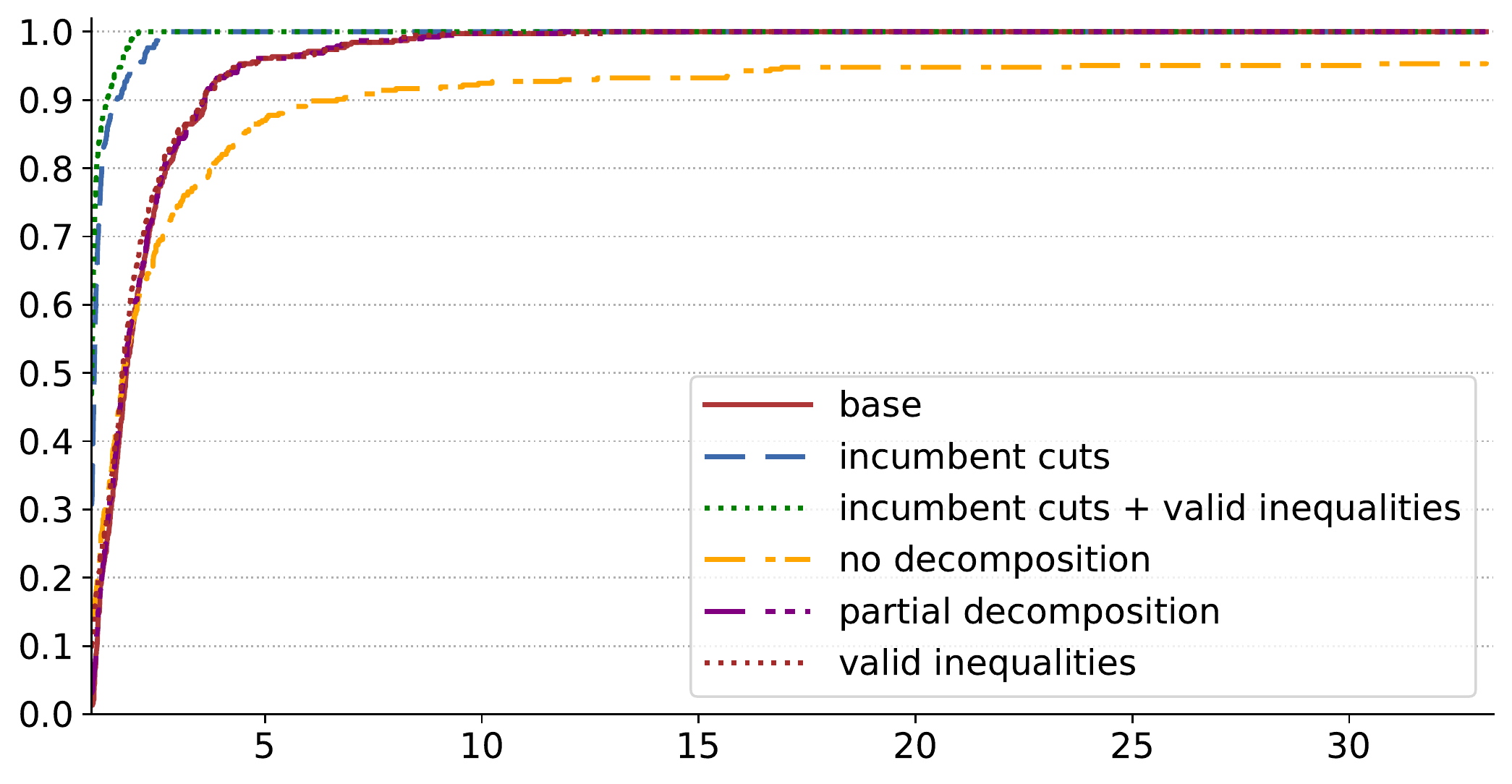}
  \caption{Performance profile: all settings on all test instances and all
    sample sizes.}
  \label{fig:perfprof-all}
\end{figure}

Figure~\ref{fig:perfprof-all} provides a comparison of the performance
profiles of all six settings on all instances and all sample sizes.
The only setting that does not always converge using the CPU budget is
\emph{no decomposition}, thus already emphasizing the need for
decomposition. We can also see that some settings are sometimes ten
times slower than others. The two settings that only generate cuts on
incumbent solutions appear to dominate the others.

\begin{figure}
  \centering
  \includegraphics[scale=.7]{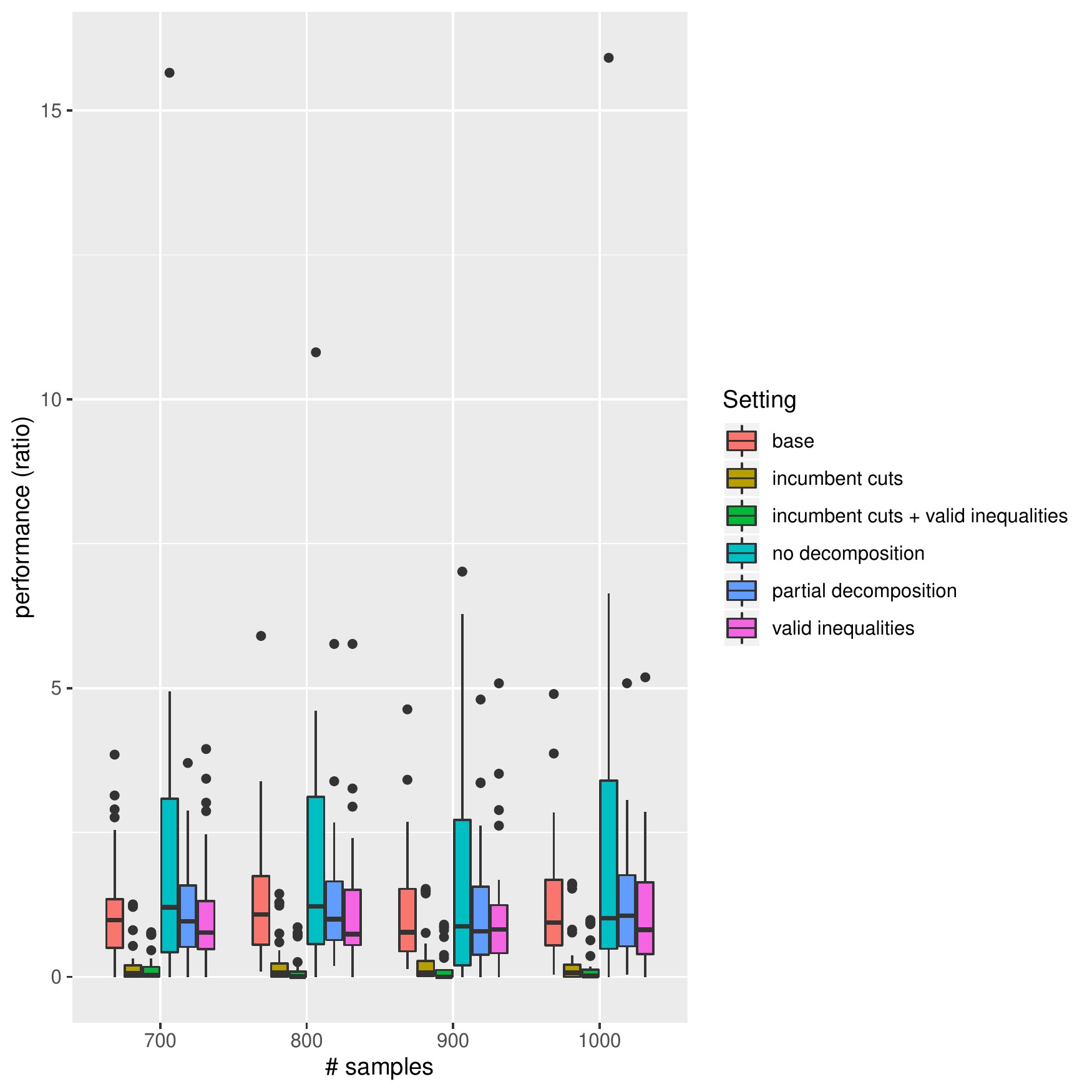}
  \caption{Box plot: all settings on all test instances, large
    sample sizes only.}
  \label{fig:boxplot-all}
\end{figure}

For further
insight, we now look at box plots for the same experimental data. We
use the \emph{ggplot2} R package~\citep{ggplot2}. Runs
for which the algorithm does not converge are discarded. For the sake
of readability, we only consider instances with at least 700 scenarios.
This box plot is depicted in Figure~\ref{fig:boxplot-all}.
We can now see that in certain cases, some methods are actually more
than 15 times slower than the best method. It appears even more
clearly that on large instances, which are the most interesting ones
since they provide a better approximation of reality, settings
generating cuts only for incumbent solutions perform better.

\begin{figure}
  \centering
  \includegraphics[scale=.7]{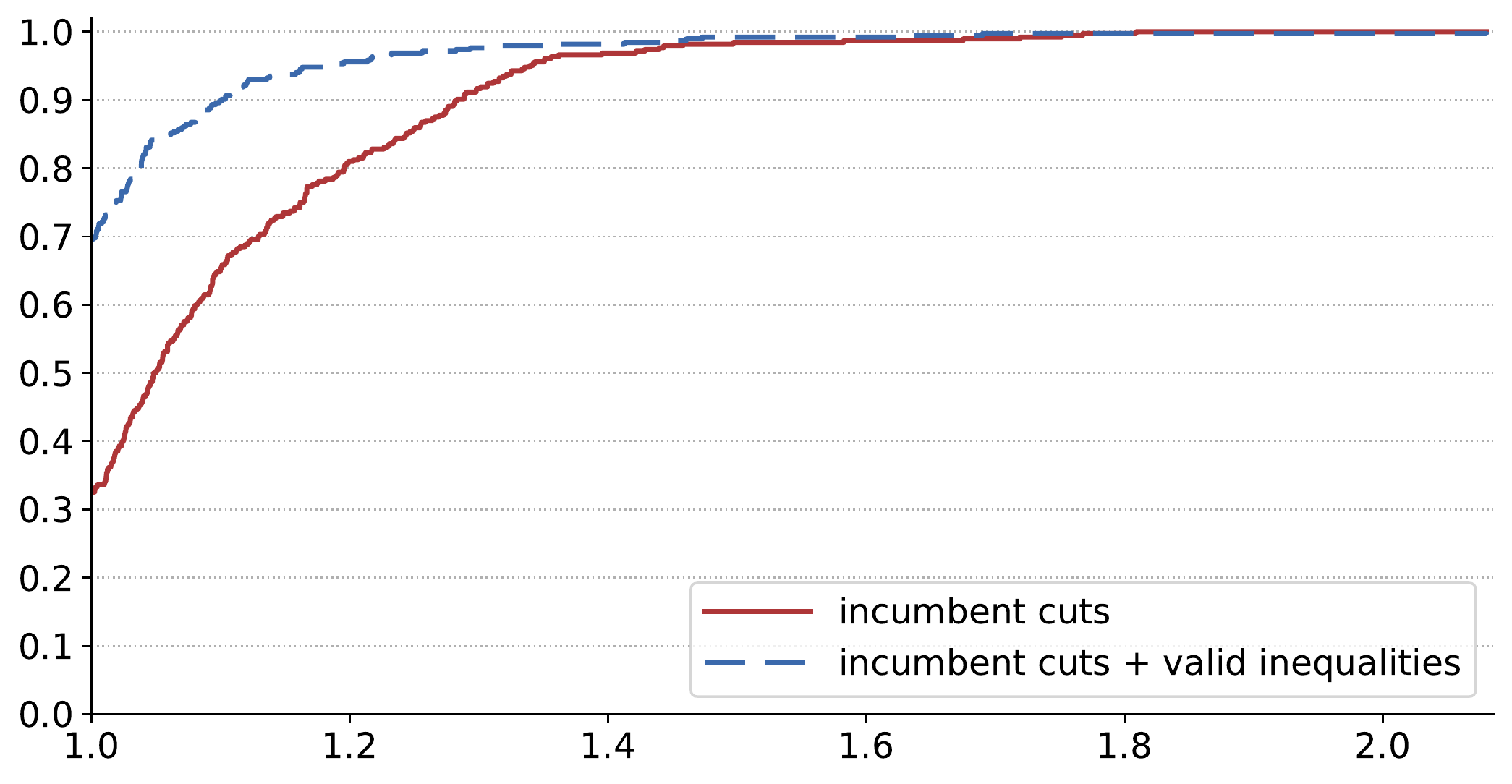}
  \caption{Performance profile: best settings on all test instances and all
    sample sizes.}
  \label{fig:perfprof-less}
\end{figure}

We now look at the two best settings only, in order to determine
whether the valid inequalities provide any kind of significant
improvement. For that purpose we first look at the performance
profiles. They are depicted in Figure~\ref{fig:perfprof-less}.
The setting that includes valid inequalities is not always
the best, as indicated by the fact that it does not start at
1. However, its curve is way above the one from the setting without
valid inequalities, indicating a better performance overall.
In general, neither setting offers any very bad performance.

\begin{figure}
  \centering
  \includegraphics[scale=.7]{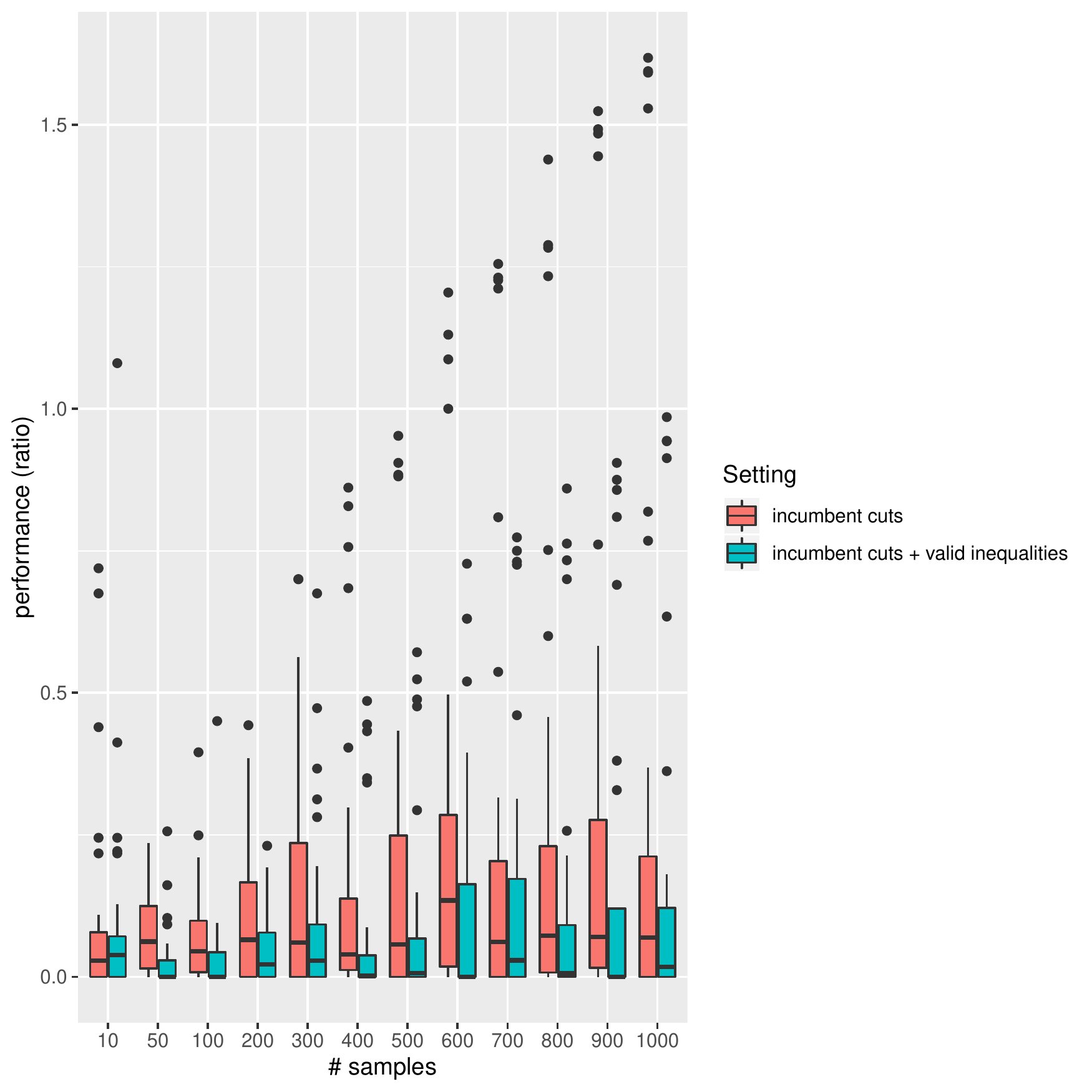}
  \caption{Box plot: best settings on all test instances and all
    sample sizes.}
  \label{fig:boxplot-less}
\end{figure}

We also provide a box plot for the two best settings in
Figure~\ref{fig:boxplot-less}.
As we can see the worst performance is below 2, meaning than no
setting is ever twice as slow as the best know setting. This, together
with prior graphics, indicates that generating cuts only at incumbent
solutions is the main cause of good performance. However, there is a
clear trend in favor of the setting that also includes valid
inequalities, observed for all sample sizes but the smallest
(10).

Based on these observations, it is clear that the best setting is the
one which both (i) only generates cuts at incumbent solutions and (ii)
includes valid inequalities in the master problem,
i.e. \emph{incumbent cuts + valid inequalities}.

\subsection{Benchmark data}
In order to facilitate future comparisons, we provide detailed results
on each instance for the overall best setting, which is
\emph{incumbent cuts + valid inequalities}. These results can be found
in Appendix~\ref{appendix:benchmark}.

\section{Conclusions and outlook}
\label{sec:conclusion}
We have defined a bi-objective facility location problem
(BOSFLP), which considers both a deterministic objective (cost
minimization) and a stochastic one (population coverage
maximization), approximated using a sampling approach. The aim of the
BOSFLP is to determine the set of efficient solutions using the Pareto
approach, but aiming for good approximations with regard to the stochastic
objective means considering large samples of random realizations, which makes
standard approaches impractical. We decomposed the original
problem and integrated Benders decomposition (L-shaped method) in a
bi-objective branch-and-bound (BIOBAB) algorithm. This is, to the best
of our knowledge, the first time that Benders decomposition has been
integrated in a bi-objective branch-and-bound algorithm. Experiments
show that the decomposition approach outperforms the explicit
consideration of all samples in the original model.

We also implemented several known improvements to the L-shaped method, and
adapted them to the bi-objective context. Among all settings, we
observed that generating Benders cuts only at the root node and at
integer solutions speeds up the search considerably. The developed 
strategy for integrating cutting plane generation into the lower bound
set algorithm generalizes to any type of cut and thus paves the way
for the development of general purpose bi-objective branch-and-cut 
algorithms relying on bound sets.
We also observed
that valid inequalities on bounds for sample-dependent values of the
stochastic objective bring a significant improvement. In both cases,
experimental observations were significant enough to justify making
these recommendations permanent, at least in the context of the
BOSFLP.

Research perspectives include the incorporation of additional
enhancements into the L-shaped method. \citet{crainic2014partial}
propose to use more sophisticated partial decomposition
strategies, such as the clustering-mean strategy in which similar
scenarios are clustered and a good representative from each cluster is
incorporated into the master program, or the convex hull strategy,
where scenarios that include other scenarios in their convex hull are
integrated in the master program.
\citet{magnanti1981accelerating} suggest improvements based on the notion of
Pareto optimal cuts, a concept which has also been successfully employed by
\citet{adulyasak2012benders}.
The literature on the single-objective L-shaped method is abundant,
and there are still lessons to be learnt in applying these techniques
to the bi-objective context. Moreover, as illustrated
by~\cite{tricoire2012bi}, in the bi-objective context there is
potential for improvements that rely on the interaction of the two
objectives; one of our future goals is to develop bi-objective
specific improvements for Benders decomposition techniques.

\paragraph{Acknowledgments}
This work has been funded by the Austrian Science Fund (FWF): P31366 and P23589. This support is gratefully acknowledged.

\bibliographystyle{apalike}
\bibliography{bosflp}

\newpage

\appendix

\section{Benchmark results}
\label{appendix:benchmark}

We report benchmark results for the best setting,
\emph{incumbent cuts + valid inequalities}. Other settings may offer a
better performance in some cases, but this is the setting that works
best overall. Each table corresponds to a different sample
size, and each table row provides indicators for one run. Indicators
for each run are the number of vertices in the instance (\emph{vertices}), the
number of linear programs solved (\emph{LPs}), the number of
branch-and-bound nodes (\emph{B\&B nodes}), the number of cuts
generated (\emph{cuts}) and the CPU effort in seconds (\emph{CPU}).

\centering
\begin{table}[h]
\begin{tabular}{lccccc}
\toprule
Instance & vertices & LPs & B\&B nodes & cuts & CPU (s) \\
\midrule
Ndiakhene & 9 & 11 & 1 & 35 & 0.19 \\
Cherif\_Lo & 10 & 143 & 21 & 78 & 0.44 \\
Thienaba & 10 & 122 & 33 & 56 & 0.35 \\
Ndieyene\_Shirak & 11 & 164 & 25 & 92 & 0.52 \\
Notto\_Gouye\_Diama & 11 & 34 & 7 & 95 & 0.3 \\
Malicounda\_Wolof & 12 & 16 & 1 & 48 & 0.24 \\
Mbayene & 12 & 220 & 41 & 69 & 0.56 \\
Pekesse & 12 & 68 & 7 & 92 & 0.4 \\
Thiadiaye & 13 & 25 & 3 & 77 & 0.31 \\
Thilmanka & 14 & 25 & 1 & 60 & 0.31 \\
Mont\_Roland & 15 & 71 & 9 & 110 & 0.54 \\
Sandira & 15 & 32 & 5 & 80 & 0.37 \\
Koul & 16 & 245 & 39 & 126 & 0.86 \\
Meouane & 17 & 137 & 31 & 125 & 0.61 \\
Neugeniene & 17 & 130 & 11 & 129 & 0.75 \\
Ndiass & 18 & 25 & 1 & 148 & 0.58 \\
Pire\_Goureye & 18 & 243 & 31 & 143 & 1.0 \\
Merina\_Dakhar & 19 & 116 & 19 & 150 & 0.86 \\
Ngandiouf & 19 & 763 & 167 & 72 & 1.71 \\
Nguekhokh & 19 & 63 & 3 & 156 & 0.77 \\
Touba\_Toul & 20 & 511 & 105 & 183 & 1.46 \\
Tassette & 21 & 425 & 79 & 252 & 1.75 \\
Diender\_Guedj & 22 & 1075 & 197 & 271 & 2.91 \\
Ndiagagniao & 24 & 4632 & 1165 & 254 & 7.26 \\
Notto & 28 & 6311 & 1591 & 262 & 10.6 \\
Pout & 29 & 5391 & 1513 & 321 & 8.91 \\
\bottomrule
\end{tabular}
\caption{Indicators for \emph{incumbent cuts + valid inequalities} (10 samples).}
\end{table}

\begin{table}[h]
\begin{tabular}{lccccc}
\toprule
Instance & vertices & LPs & B\&B nodes & cuts & CPU (s) \\
\midrule
Ndiakhene & 9 & 10 & 1 & 105 & 0.34 \\
Cherif\_Lo & 10 & 166 & 17 & 288 & 1.47 \\
Thienaba & 10 & 125 & 27 & 250 & 1.17 \\
Ndieyene\_Shirak & 11 & 443 & 83 & 279 & 2.32 \\
Notto\_Gouye\_Diama & 11 & 19 & 3 & 191 & 0.66 \\
Malicounda\_Wolof & 12 & 16 & 1 & 115 & 0.51 \\
Mbayene & 12 & 267 & 45 & 238 & 2.01 \\
Pekesse & 12 & 250 & 31 & 329 & 2.22 \\
Thiadiaye & 13 & 31 & 3 & 167 & 0.85 \\
Thilmanka & 14 & 25 & 1 & 125 & 0.75 \\
Mont\_Roland & 15 & 167 & 31 & 256 & 1.6 \\
Sandira & 15 & 67 & 17 & 179 & 1.3 \\
Koul & 16 & 205 & 31 & 296 & 2.69 \\
Meouane & 17 & 75 & 9 & 317 & 2.31 \\
Neugeniene & 17 & 115 & 11 & 362 & 2.78 \\
Ndiass & 18 & 170 & 39 & 359 & 2.46 \\
Pire\_Goureye & 18 & 322 & 67 & 369 & 4.25 \\
Merina\_Dakhar & 19 & 504 & 115 & 313 & 3.76 \\
Ngandiouf & 19 & 357 & 67 & 252 & 4.23 \\
Nguekhokh & 19 & 83 & 3 & 321 & 2.93 \\
Touba\_Toul & 20 & 968 & 203 & 541 & 7.54 \\
Tassette & 21 & 349 & 51 & 481 & 6.16 \\
Diender\_Guedj & 22 & 1131 & 231 & 678 & 9.62 \\
Ndiagagniao & 24 & 3549 & 901 & 802 & 20.21 \\
Notto & 28 & 10613 & 2733 & 928 & 55.29 \\
Pout & 29 & 4640 & 1205 & 721 & 24.47 \\
\bottomrule
\end{tabular}
\caption{Indicators for \emph{incumbent cuts + valid inequalities} (50 samples).}
\end{table}

\begin{table}[h]
\begin{tabular}{lccccc}
\toprule
Instance & vertices & LPs & B\&B nodes & cuts & CPU (s) \\
\midrule
Ndiakhene & 9 & 11 & 1 & 185 & 0.52 \\
Cherif\_Lo & 10 & 166 & 15 & 557 & 3.11 \\
Thienaba & 10 & 125 & 27 & 462 & 2.21 \\
Ndieyene\_Shirak & 11 & 306 & 63 & 601 & 4.34 \\
Notto\_Gouye\_Diama & 11 & 53 & 11 & 307 & 1.41 \\
Malicounda\_Wolof & 12 & 16 & 1 & 250 & 0.9 \\
Mbayene & 12 & 216 & 31 & 367 & 3.18 \\
Pekesse & 12 & 110 & 13 & 626 & 2.96 \\
Thiadiaye & 13 & 46 & 7 & 370 & 1.89 \\
Thilmanka & 14 & 24 & 1 & 267 & 1.42 \\
Mont\_Roland & 15 & 487 & 99 & 651 & 4.67 \\
Sandira & 15 & 37 & 5 & 311 & 2.01 \\
Koul & 16 & 228 & 37 & 678 & 5.76 \\
Meouane & 17 & 147 & 23 & 542 & 4.79 \\
Neugeniene & 17 & 232 & 31 & 597 & 6.12 \\
Ndiass & 18 & 213 & 49 & 700 & 6.2 \\
Pire\_Goureye & 18 & 412 & 93 & 783 & 8.85 \\
Merina\_Dakhar & 19 & 465 & 101 & 571 & 7.12 \\
Ngandiouf & 19 & 289 & 57 & 464 & 8.47 \\
Nguekhokh & 19 & 418 & 71 & 810 & 9.77 \\
Touba\_Toul & 20 & 2040 & 437 & 1082 & 20.34 \\
Tassette & 21 & 179 & 19 & 834 & 9.15 \\
Diender\_Guedj & 22 & 1062 & 219 & 1197 & 20.11 \\
Ndiagagniao & 24 & 3231 & 817 & 1454 & 39.72 \\
Notto & 28 & 15976 & 4073 & 1695 & 141.32 \\
Pout & 29 & 6290 & 1659 & 1486 & 62.84 \\
\bottomrule
\end{tabular}
\caption{Indicators for \emph{incumbent cuts + valid inequalities} (100 samples).}
\end{table}

\begin{table}[h]
\begin{tabular}{lccccc}
\toprule
Instance & vertices & LPs & B\&B nodes & cuts & CPU (s) \\
\midrule
Ndiakhene & 9 & 11 & 1 & 454 & 1.19 \\
Cherif\_Lo & 10 & 191 & 17 & 1127 & 7.04 \\
Thienaba & 10 & 139 & 25 & 781 & 4.48 \\
Ndieyene\_Shirak & 11 & 199 & 39 & 825 & 6.94 \\
Notto\_Gouye\_Diama & 11 & 20 & 3 & 580 & 2.37 \\
Malicounda\_Wolof & 12 & 17 & 1 & 405 & 1.88 \\
Mbayene & 12 & 252 & 43 & 745 & 7.67 \\
Pekesse & 12 & 392 & 61 & 1401 & 9.73 \\
Thiadiaye & 13 & 33 & 3 & 537 & 3.26 \\
Thilmanka & 14 & 25 & 1 & 488 & 3.34 \\
Mont\_Roland & 15 & 650 & 141 & 1471 & 12.47 \\
Sandira & 15 & 70 & 15 & 574 & 5.16 \\
Koul & 16 & 423 & 79 & 1445 & 14.95 \\
Meouane & 17 & 145 & 19 & 931 & 9.03 \\
Neugeniene & 17 & 242 & 33 & 1105 & 12.35 \\
Ndiass & 18 & 308 & 77 & 1180 & 13.19 \\
Pire\_Goureye & 18 & 720 & 155 & 1618 & 21.44 \\
Merina\_Dakhar & 19 & 560 & 127 & 1135 & 15.81 \\
Ngandiouf & 19 & 364 & 71 & 843 & 19.1 \\
Nguekhokh & 19 & 471 & 97 & 1473 & 18.5 \\
Touba\_Toul & 20 & 1419 & 275 & 2316 & 40.58 \\
Tassette & 21 & 319 & 53 & 1567 & 23.5 \\
Diender\_Guedj & 22 & 2232 & 463 & 2598 & 63.0 \\
Ndiagagniao & 24 & 4238 & 1125 & 1806 & 69.05 \\
Notto & 28 & 12136 & 3121 & 2686 & 217.35 \\
Pout & 29 & 5609 & 1469 & 2739 & 115.0 \\
\bottomrule
\end{tabular}
\caption{Indicators for \emph{incumbent cuts + valid inequalities} (200 samples).}
\end{table}

\begin{table}[h]
\begin{tabular}{lccccc}
\toprule
Instance & vertices & LPs & B\&B nodes & cuts & CPU (s) \\
\midrule
Ndiakhene & 9 & 11 & 1 & 630 & 1.62 \\
Cherif\_Lo & 10 & 181 & 23 & 1385 & 10.43 \\
Thienaba & 10 & 172 & 35 & 1281 & 5.48 \\
Ndieyene\_Shirak & 11 & 343 & 45 & 1579 & 14.72 \\
Notto\_Gouye\_Diama & 11 & 16 & 1 & 984 & 3.48 \\
Malicounda\_Wolof & 12 & 17 & 1 & 595 & 2.78 \\
Mbayene & 12 & 211 & 31 & 981 & 9.13 \\
Pekesse & 12 & 295 & 43 & 1611 & 12.07 \\
Thiadiaye & 13 & 47 & 7 & 1025 & 6.16 \\
Thilmanka & 14 & 25 & 1 & 620 & 4.98 \\
Mont\_Roland & 15 & 577 & 115 & 1683 & 17.02 \\
Sandira & 15 & 56 & 11 & 931 & 7.73 \\
Koul & 16 & 416 & 79 & 1745 & 20.66 \\
Meouane & 17 & 242 & 41 & 1419 & 16.26 \\
Neugeniene & 17 & 209 & 25 & 1573 & 20.51 \\
Ndiass & 18 & 511 & 113 & 2200 & 26.37 \\
Pire\_Goureye & 18 & 396 & 65 & 1784 & 28.09 \\
Merina\_Dakhar & 19 & 461 & 111 & 1626 & 25.53 \\
Ngandiouf & 19 & 342 & 61 & 1123 & 26.73 \\
Nguekhokh & 19 & 285 & 57 & 1868 & 23.96 \\
Touba\_Toul & 20 & 1877 & 377 & 3125 & 69.57 \\
Tassette & 21 & 255 & 33 & 2104 & 32.76 \\
Diender\_Guedj & 22 & 2077 & 417 & 3506 & 104.28 \\
Ndiagagniao & 24 & 5820 & 1409 & 4107 & 205.17 \\
Notto & 28 & 9836 & 2497 & 3624 & 302.43 \\
Pout & 29 & 6867 & 1829 & 4184 & 212.21 \\
\bottomrule
\end{tabular}
\caption{Indicators for \emph{incumbent cuts + valid inequalities} (300 samples).}
\end{table}

\begin{table}[h]
\begin{tabular}{lccccc}
\toprule
Instance & vertices & LPs & B\&B nodes & cuts & CPU (s) \\
\midrule
Ndiakhene & 9 & 11 & 1 & 720 & 2.12 \\
Cherif\_Lo & 10 & 103 & 9 & 1832 & 11.16 \\
Thienaba & 10 & 145 & 29 & 1649 & 11.08 \\
Ndieyene\_Shirak & 11 & 208 & 29 & 1968 & 14.97 \\
Notto\_Gouye\_Diama & 11 & 16 & 1 & 1055 & 4.35 \\
Malicounda\_Wolof & 12 & 17 & 1 & 926 & 4.71 \\
Mbayene & 12 & 293 & 55 & 1415 & 17.67 \\
Pekesse & 12 & 278 & 27 & 2603 & 21.68 \\
Thiadiaye & 13 & 32 & 3 & 1319 & 7.89 \\
Thilmanka & 14 & 25 & 1 & 921 & 7.47 \\
Mont\_Roland & 15 & 516 & 105 & 2108 & 23.31 \\
Sandira & 15 & 64 & 13 & 1266 & 11.55 \\
Koul & 16 & 391 & 77 & 2486 & 29.5 \\
Meouane & 17 & 208 & 35 & 1815 & 20.94 \\
Neugeniene & 17 & 165 & 23 & 2400 & 26.68 \\
Ndiass & 18 & 667 & 157 & 2731 & 39.02 \\
Pire\_Goureye & 18 & 515 & 103 & 2665 & 39.09 \\
Merina\_Dakhar & 19 & 698 & 163 & 2205 & 35.16 \\
Ngandiouf & 19 & 601 & 115 & 1633 & 45.38 \\
Nguekhokh & 19 & 420 & 109 & 2604 & 36.78 \\
Touba\_Toul & 20 & 1752 & 347 & 5580 & 110.54 \\
Tassette & 21 & 290 & 37 & 2809 & 50.14 \\
Diender\_Guedj & 22 & 1172 & 241 & 4265 & 99.05 \\
Ndiagagniao & 24 & 3502 & 835 & 4297 & 184.15 \\
Notto & 28 & 15150 & 3825 & 7627 & 724.75 \\
Pout & 29 & 7215 & 1903 & 6584 & 337.43 \\
\bottomrule
\end{tabular}
\caption{Indicators for \emph{incumbent cuts + valid inequalities} (400 samples).}
\end{table}

\begin{table}[h]
\begin{tabular}{lccccc}
\toprule
Instance & vertices & LPs & B\&B nodes & cuts & CPU (s) \\
\midrule
Ndiakhene & 9 & 11 & 1 & 870 & 2.84 \\
Cherif\_Lo & 10 & 109 & 11 & 2375 & 16.8 \\
Thienaba & 10 & 215 & 41 & 1962 & 14.68 \\
Ndieyene\_Shirak & 11 & 689 & 143 & 3782 & 35.66 \\
Notto\_Gouye\_Diama & 11 & 35 & 7 & 1325 & 7.6 \\
Malicounda\_Wolof & 12 & 17 & 1 & 999 & 5.44 \\
Mbayene & 12 & 178 & 33 & 1561 & 16.44 \\
Pekesse & 12 & 386 & 53 & 3319 & 31.52 \\
Thiadiaye & 13 & 160 & 33 & 2400 & 14.45 \\
Thilmanka & 14 & 24 & 1 & 1265 & 8.8 \\
Mont\_Roland & 15 & 632 & 131 & 3226 & 35.6 \\
Sandira & 15 & 66 & 13 & 1491 & 13.3 \\
Koul & 16 & 411 & 75 & 3569 & 49.5 \\
Meouane & 17 & 187 & 33 & 2325 & 28.24 \\
Neugeniene & 17 & 225 & 27 & 2836 & 40.15 \\
Ndiass & 18 & 465 & 139 & 2943 & 42.23 \\
Pire\_Goureye & 18 & 707 & 151 & 3290 & 53.78 \\
Merina\_Dakhar & 19 & 687 & 149 & 2490 & 41.37 \\
Ngandiouf & 19 & 746 & 159 & 1844 & 57.26 \\
Nguekhokh & 19 & 765 & 199 & 3464 & 55.49 \\
Touba\_Toul & 20 & 1446 & 293 & 5629 & 125.24 \\
Tassette & 21 & 203 & 29 & 3100 & 46.08 \\
Diender\_Guedj & 22 & 2549 & 545 & 6222 & 219.84 \\
Ndiagagniao & 24 & 3175 & 811 & 4359 & 216.4 \\
Notto & 28 & 10300 & 2627 & 7026 & 677.45 \\
Pout & 29 & 9699 & 2615 & 7750 & 531.69 \\
\bottomrule
\end{tabular}
\caption{Indicators for \emph{incumbent cuts + valid inequalities} (500 samples).}
\end{table}

\begin{table}[h]
\begin{tabular}{lccccc}
\toprule
Instance & vertices & LPs & B\&B nodes & cuts & CPU (s) \\
\midrule
Ndiakhene & 9 & 11 & 1 & 1000 & 3.22 \\
Cherif\_Lo & 10 & 194 & 21 & 2697 & 26.34 \\
Thienaba & 10 & 243 & 59 & 2480 & 15.39 \\
Ndieyene\_Shirak & 11 & 377 & 57 & 3911 & 40.69 \\
Notto\_Gouye\_Diama & 11 & 25 & 3 & 1627 & 7.89 \\
Malicounda\_Wolof & 12 & 16 & 1 & 1200 & 6.12 \\
Mbayene & 12 & 286 & 59 & 1856 & 21.41 \\
Pekesse & 12 & 117 & 11 & 3329 & 24.19 \\
Thiadiaye & 13 & 46 & 7 & 1991 & 13.96 \\
Thilmanka & 14 & 25 & 1 & 1340 & 10.53 \\
Mont\_Roland & 15 & 766 & 161 & 4389 & 56.04 \\
Sandira & 15 & 68 & 13 & 1811 & 13.83 \\
Koul & 16 & 370 & 69 & 4140 & 55.05 \\
Meouane & 17 & 269 & 47 & 2951 & 42.71 \\
Neugeniene & 17 & 170 & 17 & 3262 & 51.38 \\
Ndiass & 18 & 468 & 113 & 3462 & 53.77 \\
Pire\_Goureye & 18 & 736 & 149 & 3864 & 81.04 \\
Merina\_Dakhar & 19 & 630 & 141 & 3022 & 53.24 \\
Ngandiouf & 19 & 343 & 65 & 2549 & 68.44 \\
Nguekhokh & 19 & 893 & 225 & 4289 & 81.24 \\
Touba\_Toul & 20 & 1410 & 237 & 7485 & 185.48 \\
Tassette & 21 & 359 & 55 & 3737 & 76.48 \\
Diender\_Guedj & 22 & 1747 & 383 & 6754 & 197.8 \\
Ndiagagniao & 24 & 4169 & 1085 & 6062 & 348.91 \\
Notto & 28 & 10196 & 2545 & 7830 & 824.91 \\
Pout & 29 & 5767 & 1547 & 7242 & 417.66 \\
\bottomrule
\end{tabular}
\caption{Indicators for \emph{incumbent cuts + valid inequalities} (600 samples).}
\end{table}

\begin{table}[h]
\begin{tabular}{lccccc}
\toprule
Instance & vertices & LPs & B\&B nodes & cuts & CPU (s) \\
\midrule
Ndiakhene & 9 & 11 & 1 & 1416 & 4.94 \\
Cherif\_Lo & 10 & 199 & 21 & 3102 & 30.55 \\
Thienaba & 10 & 191 & 39 & 2310 & 22.69 \\
Ndieyene\_Shirak & 11 & 395 & 69 & 4471 & 44.67 \\
Notto\_Gouye\_Diama & 11 & 23 & 3 & 1804 & 7.97 \\
Malicounda\_Wolof & 12 & 16 & 1 & 2105 & 9.5 \\
Mbayene & 12 & 192 & 25 & 2385 & 29.63 \\
Pekesse & 12 & 241 & 35 & 4197 & 37.8 \\
Thiadiaye & 13 & 49 & 7 & 2341 & 15.59 \\
Thilmanka & 14 & 24 & 1 & 1745 & 13.53 \\
Mont\_Roland & 15 & 634 & 127 & 4459 & 55.65 \\
Sandira & 15 & 61 & 13 & 2272 & 20.31 \\
Koul & 16 & 472 & 93 & 4649 & 73.74 \\
Meouane & 17 & 148 & 23 & 3526 & 37.49 \\
Neugeniene & 17 & 338 & 59 & 3648 & 60.13 \\
Ndiass & 18 & 609 & 141 & 4299 & 82.22 \\
Pire\_Goureye & 18 & 620 & 121 & 4307 & 81.37 \\
Merina\_Dakhar & 19 & 1046 & 249 & 4454 & 88.32 \\
Ngandiouf & 19 & 418 & 73 & 2874 & 81.81 \\
Nguekhokh & 19 & 1112 & 281 & 5363 & 94.39 \\
Touba\_Toul & 20 & 3078 & 677 & 7367 & 254.06 \\
Tassette & 21 & 365 & 61 & 4757 & 90.2 \\
Diender\_Guedj & 22 & 1862 & 393 & 7655 & 271.29 \\
Ndiagagniao & 24 & 3690 & 937 & 6732 & 362.6 \\
Notto & 28 & 12753 & 3273 & 9053 & 1153.8 \\
Pout & 29 & 8741 & 2339 & 12372 & 892.39 \\
\bottomrule
\end{tabular}
\caption{Indicators for \emph{incumbent cuts + valid inequalities} (700 samples).}
\end{table}

\begin{table}[h]
\begin{tabular}{lccccc}
\toprule
Instance & vertices & LPs & B\&B nodes & cuts & CPU (s) \\
\midrule
Ndiakhene & 9 & 11 & 1 & 1390 & 4.67 \\
Cherif\_Lo & 10 & 214 & 23 & 4416 & 37.43 \\
Thienaba & 10 & 182 & 37 & 2525 & 20.82 \\
Ndieyene\_Shirak & 11 & 253 & 43 & 4035 & 36.19 \\
Notto\_Gouye\_Diama & 11 & 45 & 9 & 2249 & 12.34 \\
Malicounda\_Wolof & 12 & 17 & 1 & 1615 & 8.8 \\
Mbayene & 12 & 219 & 37 & 2676 & 35.43 \\
Pekesse & 12 & 103 & 13 & 3857 & 29.22 \\
Thiadiaye & 13 & 49 & 7 & 2735 & 19.43 \\
Thilmanka & 14 & 24 & 1 & 1580 & 14.12 \\
Mont\_Roland & 15 & 582 & 113 & 4445 & 58.15 \\
Sandira & 15 & 83 & 17 & 2110 & 24.33 \\
Koul & 16 & 472 & 91 & 5776 & 93.74 \\
Meouane & 17 & 237 & 39 & 3863 & 58.61 \\
Neugeniene & 17 & 258 & 39 & 3767 & 75.57 \\
Ndiass & 18 & 330 & 85 & 4427 & 66.73 \\
Pire\_Goureye & 18 & 554 & 111 & 4735 & 80.23 \\
Merina\_Dakhar & 19 & 357 & 73 & 3849 & 63.43 \\
Ngandiouf & 19 & 572 & 107 & 3332 & 103.31 \\
Nguekhokh & 19 & 784 & 223 & 5581 & 93.88 \\
Touba\_Toul & 20 & 2121 & 427 & 8548 & 258.43 \\
Tassette & 21 & 297 & 39 & 5746 & 110.27 \\
Diender\_Guedj & 22 & 2298 & 487 & 8687 & 351.12 \\
Ndiagagniao & 24 & 4547 & 1109 & 9431 & 548.11 \\
Notto & 28 & 14639 & 3699 & 11121 & 1522.24 \\
Pout & 29 & 5104 & 1275 & 10398 & 651.21 \\
\bottomrule
\end{tabular}
\caption{Indicators for \emph{incumbent cuts + valid inequalities} (800 samples).}
\end{table}

\begin{table}[h]
\begin{tabular}{lccccc}
\toprule
Instance & vertices & LPs & B\&B nodes & cuts & CPU (s) \\
\midrule
Ndiakhene & 9 & 11 & 1 & 1555 & 5.21 \\
Cherif\_Lo & 10 & 117 & 11 & 3715 & 34.87 \\
Thienaba & 10 & 228 & 49 & 3486 & 30.81 \\
Ndieyene\_Shirak & 11 & 385 & 67 & 4937 & 51.89 \\
Notto\_Gouye\_Diama & 11 & 63 & 13 & 2352 & 14.58 \\
Malicounda\_Wolof & 12 & 17 & 1 & 2135 & 10.49 \\
Mbayene & 12 & 295 & 59 & 3287 & 42.85 \\
Pekesse & 12 & 315 & 45 & 4949 & 54.64 \\
Thiadiaye & 13 & 43 & 7 & 3120 & 22.09 \\
Thilmanka & 14 & 25 & 1 & 2040 & 17.46 \\
Mont\_Roland & 15 & 906 & 201 & 7784 & 113.98 \\
Sandira & 15 & 77 & 17 & 2643 & 27.51 \\
Koul & 16 & 448 & 87 & 5799 & 102.2 \\
Meouane & 17 & 292 & 53 & 4488 & 76.73 \\
Neugeniene & 17 & 245 & 35 & 5181 & 85.7 \\
Ndiass & 18 & 554 & 129 & 6010 & 108.32 \\
Pire\_Goureye & 18 & 446 & 81 & 6068 & 123.23 \\
Merina\_Dakhar & 19 & 561 & 113 & 4709 & 107.37 \\
Ngandiouf & 19 & 474 & 87 & 3664 & 107.18 \\
Nguekhokh & 19 & 520 & 151 & 5168 & 89.1 \\
Touba\_Toul & 20 & 2233 & 449 & 12664 & 380.63 \\
Tassette & 21 & 489 & 85 & 5656 & 134.4 \\
Diender\_Guedj & 22 & 2848 & 617 & 11535 & 492.3 \\
Ndiagagniao & 24 & 4065 & 1009 & 10882 & 642.0 \\
Notto & 28 & 9397 & 2361 & 11463 & 1171.44 \\
Pout & 29 & 7902 & 2129 & 14670 & 1056.83 \\
\bottomrule
\end{tabular}
\caption{Indicators for \emph{incumbent cuts + valid inequalities} (900 samples).}
\end{table}

\begin{table}[h]
\begin{tabular}{lccccc}
\toprule
Instance & vertices & LPs & B\&B nodes & cuts & CPU (s) \\
\midrule
Ndiakhene & 9 & 11 & 1 & 1725 & 5.56 \\
Cherif\_Lo & 10 & 123 & 17 & 5043 & 40.84 \\
Thienaba & 10 & 172 & 35 & 3807 & 34.49 \\
Ndieyene\_Shirak & 11 & 288 & 51 & 6377 & 62.43 \\
Notto\_Gouye\_Diama & 11 & 20 & 3 & 2491 & 12.72 \\
Malicounda\_Wolof & 12 & 16 & 1 & 2190 & 11.7 \\
Mbayene & 12 & 305 & 59 & 3728 & 39.78 \\
Pekesse & 12 & 154 & 15 & 5669 & 46.5 \\
Thiadiaye & 13 & 268 & 59 & 5460 & 40.66 \\
Thilmanka & 14 & 25 & 1 & 2445 & 20.35 \\
Mont\_Roland & 15 & 808 & 169 & 7302 & 106.04 \\
Sandira & 15 & 59 & 13 & 2846 & 29.49 \\
Koul & 16 & 408 & 79 & 6350 & 108.44 \\
Meouane & 17 & 277 & 51 & 5350 & 83.05 \\
Neugeniene & 17 & 182 & 23 & 4886 & 81.65 \\
Ndiass & 18 & 332 & 87 & 5734 & 88.5 \\
Pire\_Goureye & 18 & 559 & 113 & 6854 & 138.85 \\
Merina\_Dakhar & 19 & 448 & 91 & 4683 & 85.16 \\
Ngandiouf & 19 & 660 & 139 & 4171 & 129.67 \\
Nguekhokh & 19 & 717 & 199 & 6110 & 121.03 \\
Touba\_Toul & 20 & 1845 & 397 & 8545 & 292.24 \\
Tassette & 21 & 305 & 51 & 7060 & 154.45 \\
Diender\_Guedj & 22 & 1714 & 385 & 10891 & 389.12 \\
Ndiagagniao & 24 & 3906 & 915 & 11495 & 793.63 \\
Notto & 28 & 13451 & 3429 & 13171 & 2019.97 \\
Pout & 29 & 7246 & 1965 & 15670 & 1159.21 \\
\bottomrule
\end{tabular}
\caption{Indicators for \emph{incumbent cuts + valid inequalities} (1000 samples).}
\end{table}

\end{document}